\documentclass[12pt,reqno,a4paper]{amsart} 
 
\usepackage{mathtools}
\usepackage[breaklinks]{hyperref} 
\usepackage{verbatim}
\usepackage[headings]{fullpage} 

\usepackage{microtype}  
 
\renewcommand{\qedsymbol}{$\square$}
\newenvironment{Proof}[1][Proof]{\par\noindent\textbf{#1.}~}
{\hfill\qedsymbol\smallskip\par}

\newcommand{\dx}{\mathrm{d}}
\newcommand{\eps}{\varepsilon}

\newcommand{\Odip}[2]{\mathcal{O}_{#1}\!\left(#2\right)\mathchoice{\!}{}{}{}}

\newcommand{\Odipm}[2]{\mathcal{O}_{#1}\bigl(#2\bigr)\mathchoice{\!}{}{}{}}
\newcommand{\Odig}[1]{\mathcal{O}\Bigl(#1\Bigr)\mathchoice{\!}{}{}{}}

\newcommand{\Odim}[1]{\mathcal{O}\bigl(#1\bigr)}
\newcommand{\Odi}[1]{\Odip{}{#1}}
\newcommand{\odip}[2]{{o}_{#1}\!\left(#2\right)\mathchoice{\!}{}{}{}}
\newcommand{\odi}[1]{\odip{}{#1}}
 
\makeatletter

\def\env@Biggcases{%
  \let\@ifnextchar\new@ifnextchar
  \Biggl\lbrace
  \def\arraystretch{1.2}%
  \array{@{}l@{\quad}l@{}}%
}
\makeatother

\newtheoremstyle{sltheorems}
{10pt}
{6pt}
{\slshape}
{}
{\bfseries}
{.}
{.5em}
{\thmname{#1}\thmnumber{ #2}\thmnote{ (#3)}}

\theoremstyle{sltheorems} 
\newtheorem{Theorem}{Theorem}
\newtheorem{Lemma}{Lemma}

\allowdisplaybreaks
\title[Goldbach average in short intervals]{Ces\`aro average in short intervals for Goldbach numbers}
\author[Languasco \lowercase{and}  Zaccagnini]{Alessandro Languasco \lowercase{and} Alessandro Zaccagnini}


\begin{document}
%
\begin{abstract}
 Let $\Lambda$ be the von Mangoldt function and 
\(
R(n) = \sum_{h+k=n} \Lambda(h)\Lambda(k).
\)
Let further $N,H$ be two integers, $N\ge 2$, $1\le H \le N$, and 
assume  that the Riemann Hypothesis  holds.
Then
\begin{align*}
\sum_{n=N-H}^{N+H}
R(n) \Bigl(1- \frac{\vert n- N \vert}{H}\Bigr)
&=
HN
-\frac{2}{H}
\sum_{\rho} \frac{(N+H)^{\rho +2} - 2 N^{\rho+2} +(N-H)^{\rho +2} }{\rho (\rho + 1)(\rho + 2)}
\\&
+
\Odig{N \Bigl(\log \frac{2N}{H}\Bigr)^2 + H (\log N)^2 \log(2H) 
}\ ,
\end{align*}
where $\rho=1/2+i\gamma$ runs over
the non-trivial zeros of the Riemann zeta function $\zeta(s)$.

\end{abstract}
\maketitle
\section{Introduction}

Let $\Lambda$ be the von Mangoldt function and
\[
R(n) = \sum_{h_1+h_2=n} \Lambda(h_1)\Lambda(h_2)
\]
be the counting function for the Goldbach numbers.
In this paper we are looking for an explicit formula 
for a Ces\`aro average of $R(n)$ in short intervals.
Concerning long intervals, we should mention 
our result in \cite{LanguascoZ2012a}: assuming 
the Riemann Hypothesis (RH)  
we have
\[ 
 \sum_{n \le N} 
R(n) 
=
\frac{N^{2}}{2}
-2
\sum_{\rho} \frac{N^{\rho + 1}}{\rho (\rho + 1)}
+
\Odim{N (\log N)^{3}
},
\]
where  $N $ is a  large integer and $\rho=1/2+i\gamma$ runs over
the non-trivial zeros of the Riemann zeta function $\zeta(s)$.
We also mention its extension to the Ces\`aro average case
by Goldston-Yang \cite{GoldstonY2016} again under the assumption of RH:
\begin{equation}
\label{Goldston-Yang}
   \sum_{n \le N}  R(n) \Bigl(1 - \frac{n}{N}\Bigr)
  =
  \frac{N^{2}}{6}
  -
  2
  \sum_{\rho} \frac{N^{\rho+1}}{\rho(\rho+1)(\rho+2)}   +  \Odim{N} .
\end{equation}
We also recall our unconditional result in \cite{LanguascoZ2015a}, see
also \cite{Languasco2016a}:  let $k>1$ be a real number; we have
\begin{align*}
\sum_{n \le N} R(n) \frac{(1 - n/N)^k}{\Gamma(k + 1)}
  &=
  \frac{N^{2}}{\Gamma(k + 3)}
  -
  2
  \sum_{\rho} \frac{\Gamma(\rho)}{\Gamma(\rho + k + 2)} N^{\rho+1} \\
  &\qquad+
\label{expl-form-Goldbach}
  \sum_{\rho_1} \sum_{\rho_2}
    \frac{\Gamma(\rho_1) \Gamma(\rho_2)}{\Gamma(\rho_1 + \rho_2 + k + 1)}
    N^{\rho_1 + \rho_2}
  +
  \Odipm{k}{N}\ ,
  \end{align*}
where $\rho_1,\rho_2$  run  over
the non-trivial zeros of the Riemann zeta function $\zeta(s)$
and $\Gamma(s)$ is Euler's function.
Our result here is
\begin{Theorem}
\label{Main-Th}
Let $N,H$ be two integers, $N\ge 2$, $1\le H \le N$. Assume  that the Riemann Hypothesis (RH) holds.
Then
\begin{align*}
\sum_{n=N-H}^{N+H}
R(n) \Bigl(1- \frac{\vert n- N \vert}{H}\Bigr)
&=
HN
-\frac{2}{H}
\sum_{\rho} \frac{(N+H)^{\rho +2} - 2 N^{\rho+2} +(N-H)^{\rho +2} }{\rho (\rho + 1)(\rho + 2)}
\\&
+
\Odig{  N \Bigl(\log \frac{2N}{H}\Bigr)^2 + H (\log N)^2 \log(2H) 
}\ ,
\end{align*}
where $\rho=1/2+i\gamma$ runs over
the non-trivial zeros of the Riemann zeta function $\zeta(s)$.
\end{Theorem}
 
 The second difference involved in the zero-depending term is natural
 since it depends on the symmetric nature of the short-interval
 Ces\`aro weight used in Theorem \ref{Main-Th}. Its unconditional 
 order of magnitude is $\ll H N \exp(-c_1 (\log N)^{3/5}(\log \log n)^{-1/5}) +N$,
where $c_1>0$ is an absolute constant, while,
under the assumption of RH, it is $\ll HN^{1/2} (\log N)^2+N$,
 see Section \ref{S-order}.

In fact  we will obtain Theorem \ref{Main-Th} as a consequence
of a weighted result. Letting $\psi(x)=\sum_{m\leq x} \Lambda(m)$, we have
\begin{Theorem}
\label{average-Th}
Let $N,H$ be two integers, $N\ge 2$, $1\le H \le N$ and $y \in [-H, H]$. Assume that 
the Riemann Hypothesis (RH) holds.
Then
\begin{equation}
\label{average}
\max_{y \in [-H, H)}\ 
\Bigl\vert
\sum_{n=N-H}^{N+y} 
e^{-n/N}
\Bigl(
R(n) - (2\psi(n)-n)
\Bigr)
\Bigl(1- \frac{\vert n- N \vert}{H}\Bigr)
\Bigr\vert
\ll
 N  (\log N )^2  \log (2H)
\end{equation}
and
\begin{equation}
\label{average-H} 
\Bigl\vert
\sum_{n=N-H}^{N+H} e^{-n/N}
\Bigl(
R(n) - (2\psi(n)-n)
\Bigr)
\Bigl(1- \frac{\vert n- N \vert}{H}\Bigr)
\Bigr\vert
\ll
 N   \Bigl(\log \frac{2N}{H}\Bigr)^2.
\end{equation}
\end{Theorem}
%
The better estimate for  the case $y=H$ depends on
the second point of Lemma \ref{T-behaviour} below
in which we have a more efficient estimate 
for the exponential sum $T_H(H,H;\alpha)$, defined in \eqref{T-def},
attached to the Ces\`aro weight.

For $H=N$ we can compare Theorem \ref{Main-Th} with \eqref{Goldston-Yang}
and it is clear that the previously mentioned 
weakness of the available estimates for $T_H(H,y;\alpha)$, again defined in \eqref{T-def},
when $y\neq H$ leads us to a weaker final estimate by a factor $(\log N)^3$. 
Unfortunately it seems that Lemma \ref{T-behaviour} is optimal, see 
the remark after its proof, and hence this is a serious limitation for our method.

After being shown this paper, Goldston \& Yang told us that it should
be possible to combine their technique in \cite{GoldstonY2016} with
our Lemmas \ref{sum-integral-modificato} and \ref{sum-integral-pesato}
below to remove the second error term in the statement of Theorem
\ref{Main-Th}.

In order to match the case $H=N$ with our method,
we should have a  more efficient way of removing the $e^{-n/N}$ weight
(which naturally arises from the use of infinite series, see \eqref{tildeS-def});
unfortunately the partial summation strategy we used to achieve this goal
needs a uniform result on $y$. This leads to the first estimate in Lemma
\ref{T-behaviour} and hence  our global method 
is efficient essentially only for $H\ll N (\log \log N)/(\log N)^3$.

As we did in \cite{LanguascoZ2012a},
 we will  use the original Hardy and Littlewood \cite{HardyL1923} circle
method setting, \emph{i.e.},
the weighted exponential sum
\begin{equation}
\label{tildeS-def}
\widetilde{S}(\alpha)
=
\sum_{n=1}^{\infty}
\Lambda(n) e^{-n/N}
e(n\alpha),
\end{equation}
where $e(x)=\exp(2\pi i x)$.
Such a function was also used by Linnik \cite{Linnik1946,Linnik1952}.

\medskip
\textbf {Acknowledgments.} 
We thank the referee for pointing out several inaccuracies in a previous version of this paper.

\section{Setting of the circle method}

For brevity, throughout the paper we write
\begin{equation}
\label{def-z}
  z
  =
  \frac1N - 2 \pi i \alpha,
\end{equation}
where $N$ is a large integer and $\alpha \in [-1/2, 1/2]$.
The first lemma is an $L^{2}$-estimate for the
difference  $\widetilde{S}(\alpha) - 1/z$.
\begin{Lemma} 
\label{LP-Lemma-mod}
Assume RH. Let $N$ be a
sufficiently large integer and $z$ be as in \eqref{def-z}.
For  $0 \leq \xi \leq 1/2$, we have
\[
\int_{-\xi}^{\xi}
\Bigl\vert
\widetilde{S}(\alpha) - \frac{1}{z}
\Bigr\vert^{2}
\dx \alpha
\ll
N\xi \bigl(1+ \log (2N\xi)\bigr)^{2}.
\]
\end{Lemma}

\begin{Proof}
This follows immediately from the proof of
Theorem 1 of \cite{LanguascoP1994}.
We just have to pay attention to the final estimate 
of eq. (22) on page 315 there. A slightly more
careful estimate immediately gives that (22) 
can be replaced by
\[
\ll \sum_{\gamma_1>0}
\exp\Bigl(-\frac{c}{2}\frac{\gamma_1}{N\eta}\Bigr)
\sum_{\gamma_2>0}\frac{1}{1+\vert \gamma_1-\gamma_2\vert ^2}
 \ll 
 N\eta \bigl( \log (2N\eta)\bigr)^{2}.
\]
The final estimate follows at once.
\end{Proof}

The next four lemmas do not depend on RH.
By the residue theorem one can obtain
\begin{Lemma}[Eq.~(29) of \cite{LanguascoP1994}]
\label{residue}
Let $N \geq 2$ and $1\leq n \leq 2N$ be integers; let further $z$ be as in \eqref{def-z}.
We have
\[
\int_{-\frac{1}{2}}^{\frac{1}{2}}
\frac{e(-n\alpha)}{z^2}
\ \dx \alpha
=
ne^{-n/N} + \Odi{1}
\]
uniformly for every $n \leq 2N$.
\end{Lemma}

\begin{Lemma}[Lemma 2.3  of \cite{LanguascoZ2012a}]
\label{incond-mean-square}
Let $N$ be a sufficiently large integer and $z$ be as in \eqref{def-z}.
We have
\[
\int_{-\frac{1}{2}}^{\frac{1}{2}}
\Bigl\vert
\widetilde{S}(\alpha) - \frac{1}{z}
\Bigr\vert^{2}
\ \dx \alpha
=
\frac{N}{2} \log N
+
\Odig{N (\log N)^{1/2}}.
\]
\end{Lemma}
   
Let
\begin{equation}
\label{V-def}
  V(\alpha)
  =
  \sum_{m = 1}^{\infty} e^{-m / N} e(m \alpha)
  =
  \sum_{m = 1}^{\infty} e^{-m z}
  =
  \frac1{e^z - 1}.
\end{equation}

\begin{Lemma}[Lemma 2.4 of \cite{LanguascoZ2012a}]
\label{V-behaviour}
If $z$ satisfies \eqref{def-z} then $V(\alpha) = z^{-1} + \Odi{1}$.
\end{Lemma}

Let now 
\begin{equation}
\label{T-def}
t_H(m)=H-\vert m \vert 
\quad
\textrm{and}
\quad
T_H(N,y;\alpha) =
\sum_{n=N-H}^{N+y} t_H(n-N) e(n\alpha).
\end{equation}

\begin{Lemma} 
\label{T-behaviour}
Let $N,H$ be two integers, $N\ge 2$, $1\le H \le N$.
For every $y\in [-H,H)$ and $\alpha\in[-1/2,1/2]$,  we have  
\[
T_H(N,y;\alpha)  
\ll H \min\Bigl(H ; \frac{1}{\Vert \alpha\Vert}\Bigr).
\]
Moreover,  for every   $\alpha\in[-1/2,1/2]$, we also have
\[
T_H(N,H;\alpha)
\ll \min\Bigl(H^2; \frac{1}{\Vert \alpha\Vert^2}\Bigr).
\]
\end{Lemma}
\begin{Proof}
First of all we recall the well-known estimate
\begin{equation}
\label{expsum-semplice}
\sum_{m=1}^{u}  e(m\alpha) \ll \min\Bigl(u ; \frac{1}{\Vert \alpha\Vert}\Bigr).
\end{equation}
Let now $y\in [-H,H)$. Then 
\begin{equation}
\label{first-case-estim1}
\vert T_H(N,y;\alpha)  \vert \le 
\sum_{n=N-H}^{N+y} t_H(n-N) 
\ll H(H+y+1) \ll H^2.
\end{equation}
Moreover if $y \ge 0$ we get 
\begin{equation}
\label{ABC-def}
T_H(N,y;\alpha) = 
H \sum_{n=N-H}^{N+y}  e(n\alpha) 
- 
  \sum_{m=0}^{y} m e((N+m)\alpha) 
-
 \sum_{m=1}^{H} m e((N-m)\alpha)
=A- B - C,
\end{equation}
say.
By partial summation and \eqref{expsum-semplice} we get
\begin{equation}
\label{first-case-estim2}
B=  y \sum_{m=0}^{y}   e((N+m)\alpha)  - \int _0^y \sum_{m=0}^w   e((N+m)\alpha)\ \dx w
\ll 
\frac{ y  }{\Vert \alpha\Vert} + 
 \int _0^y \frac{\dx w}{\Vert \alpha\Vert}  
 \ll
 \frac{ y  }{\Vert \alpha\Vert}
 \ll
 \frac{ H  }{\Vert \alpha\Vert}.
\end{equation}
Arguing analogously we have 
\begin{equation}
\label{first-case-estim3}
 C \ll \frac{ H  }{\Vert \alpha\Vert},
 \end{equation}
while the inequality $A \ll  H  /\Vert \alpha\Vert$ follows from \eqref{expsum-semplice}.
If $y<0$ then we can write that
\[
T_H(N,y;\alpha) = 
H \sum_{n=N-H}^{N+y}  e(n\alpha) 
-
 \sum_{m=-y}^{H} m e((N-m)\alpha)
=A- D ,
\]
say,
where $A$ is defined in \eqref{ABC-def}.
Arguing as we did for $B$ we get
\begin{equation}
\label{first-case-estim4}
D \ll \frac{H}{\Vert \alpha\Vert}.
\end{equation}
Combining \eqref{first-case-estim1}-\eqref{first-case-estim4}
the first part of the lemma follows for every $y\in [-H,H)$.
The second part of the lemma follows by \eqref{expsum-semplice} 
and the fact that in this case we can write
\[
T_H(N,H;\alpha)
= 
\sum_{m=-H}^{H}
t_H(m) e(m\alpha) e(N\alpha) 
= 
\Bigl\vert \sum_{m=1}^{H}  e(m\alpha)\Bigr\vert ^2  e(N\alpha).
\eqno{}\qed
\]
\renewcommand{\qedsymbol}{}
\end{Proof}
\textbf{Remark.}
We remark that the estimate for $T_H(N,y;\alpha)$, $y\ne H$, is
essentially optimal. 
For brevity, we only deal with $T_H(N; y, \alpha)$ for $y \in [0, H]$.
It is not hard to prove by induction that
\begin{align*}
  T_H(N; y, \alpha)
  &=
  e(N \alpha)
  \sum_{n = -H}^y (H - |n|) e(n \alpha) \\
  &=
  \frac{e( (N + y + 1) \alpha)}{1 - e(\alpha)}
  \cdot
  (y - H)
  +
  \frac{e( (N + 1) \alpha)}{(1 - e(\alpha))^2}
  \cdot
  \bigl( e(y \alpha) - 2 + e(- H \alpha) \bigr).
\end{align*}
In the critical range $\alpha \in [H^{-1}, 1 / 2]$ the last summand
has a smaller order of magnitude than $H\Vert \alpha \Vert^{-1}$, 
and this implies that the
bound $T_H(N; y, \alpha) \ll H \Vert \alpha \Vert^{-1}$ is sharp, at
least when $y \le H / 2$, say.

\medskip
We build now the zero-depending term we have in Theorem \ref{Main-Th}.
The first step is the following
\begin{Lemma}
\label{I2-lemma}
Let $N,H$ be two integers, $N\ge 2$, $1\le H \le N$ and $z$ be as in
\eqref{def-z}.
For every $y\in [-H,H)$
we have
\begin{align}
\notag
  \int_{-1/2}^{1/2} T_H(N,y;-\alpha)
\frac{(\widetilde{S}(\alpha)-1/z)}z \, \dx \alpha
&=
\sum_{n=N-H}^{N+y} e^{-n / N} t_H(n-N)  (\psi(n) - n)\\
\label{I2-eval}
&
  +
  \Odi{H^{3/2} N^{1/2} (\log N)^{1/2}}.
\end{align}
\end{Lemma}
We remark that Lemma \ref{I2-lemma}, which is a modified version
of Lemma 2.5 of \cite{LanguascoZ2012a}, is unconditional
and hence it implies, using also Lemmas \ref{sum-integral-modificato}-\ref{sum-integral-pesato}, 
that the ability of detecting the zero-depending term
of the Riemann zeta function  $\zeta(s)$
in Theorem \ref{Main-Th} does not depend on RH.

\begin{Proof}
Writing
$\widetilde{R}(\alpha)=\widetilde{S}(\alpha) - 1/z$,
by Lemma \ref{V-behaviour} we have
\begin{align} 
\notag
  \int_{-1/2}^{1/2}\!\!\!   T_H(N,y;-\alpha)& \frac{\widetilde{R}(\alpha)}z \, \dx \alpha
   =
  \int_{-1/2}^{1/2} \!\!\! T_H(N,y;-\alpha)\widetilde{R}(\alpha) V(\alpha) \, \dx \alpha
+
  \Odig{\int_{-1/2}^{1/2} \vert  T_H(N,y;-\alpha)\vert  \, \vert \widetilde{R}(\alpha)\vert  \, \dx \alpha}
  \\
    & =
  \int_{-1/2}^{1/2} T_H(N,y;-\alpha)\widetilde{R}(\alpha) V(\alpha) \, \dx \alpha
\label{first-step}  +
  \Odi{  (H^3 N \log N)^{1/2}},
\end{align}
since, by Lemmas \ref{incond-mean-square} and \ref{T-behaviour},
the error term above is
\[
  \ll
  \Bigl( \int_{-1/2}^{1/2} \vert T_H(N,y;-\alpha)\vert ^2 \, \dx \alpha \Bigr)^{1/2}
  \Bigl( \int_{-1/2}^{1/2} \vert \widetilde{R}(\alpha)\vert ^2 \, \dx \alpha \Bigr)^{1/2}
  \ll
  (H^3 N \log N)^{1/2}.
\]
Again by Lemma \ref{V-behaviour}, we have
\[
  \widetilde{R}(\alpha)
  =
  \widetilde{S}(\alpha)
  -
  \frac1z
  =
   \widetilde{S}(\alpha)
  -
  V(\alpha)
  +
  \Odi{1}
\]
and hence \eqref{first-step} implies
\begin{align}
\notag
  \int_{-1/2}^{1/2}  T_H(N,y;-\alpha)& \frac{\widetilde{R}(\alpha)}z \, \dx \alpha
 =
  \int_{-1/2}^{1/2}
   T_H(N,y;-\alpha)
    \Bigl(\widetilde{S}(\alpha) - V(\alpha) \bigr) V(\alpha) \, \dx \alpha \\
\label{second-step}
  &
\quad
+
  \Odig{\int_{-1/2}^{1/2} \vert T_H(N,y;-\alpha)\vert  \,  \vert V(\alpha)\vert  \, \dx \alpha}
  +
  \Odi{  (H^3 N \log N)^{1/2}}.
\end{align}
The Cauchy-Schwarz inequality, Lemma \ref{T-behaviour} and the Parseval theorem imply that
\begin{align}
\notag
  \int_{-1/2}^{1/2} \vert T_H(N,y;-\alpha)\vert  \, \vert V(\alpha)\vert  \, \dx \alpha
  &\le
  \Bigl( \int_{-1/2}^{1/2} \vert T_H(N,y;-\alpha))\vert ^2 \, \dx \alpha \Bigr)^{1/2}
  \Bigl( \int_{-1/2}^{1/2} \vert V(\alpha)\vert ^2 \, \dx \alpha \Bigr)^{1/2} \\
  &\ll
\label{err-second-step}
  \Bigl( H^3 \sum_{m = 1}^{\infty} e^{-2 m / N} \Bigr)^{1/2}
  \ll
  (H^3 N)^{1/2}.
\end{align}
By \eqref{second-step}-\eqref{err-second-step}, we have
\begin{align}
\notag
  \int_{-1/2}^{1/2}  T_H(N,y;-\alpha)\frac{\widetilde{R}(\alpha)}z \, \dx \alpha
  &=
  \int_{-1/2}^{1/2}
    T_H(N,y;-\alpha)
    \Bigl(\widetilde{S}(\alpha) - V(\alpha) \bigr) V(\alpha) \, \dx \alpha
    \\
    \label{start-third-step}
    &
  +
  \Odi{(H^3 N \log N)^{1/2}}.
\end{align}
Now, by \eqref{tildeS-def} and \eqref{V-def}, we can write
\[
  \widetilde{S}(\alpha)
  -
  V(\alpha)
  =
  \sum_{m = 1}^{\infty} (\Lambda(m) - 1) e^{-m / N}  e(m \alpha)
\]
so that
\begin{align}
\notag
   \int_{-1/2}^{1/2}&
    T_H(N,y;-\alpha)
    \Bigl(\widetilde{S}(\alpha) - V(\alpha) \bigr) V(\alpha) \, \dx \alpha \\
\notag
  &=
  \sum_{m= -H}^{y} t_H(m)
 \sum_{m_1 = 1}^{\infty} (\Lambda(m_1) - 1) e^{-m_1 / N}
   \sum_{m_2 = 1}^{\infty} e^{-m_2 / N}
   \int_{-1/2}^{1/2} e((m_1 + m_2 - m -N) \alpha) \, \dx \alpha \\
\notag
 &=
  \sum_{m= -H}^{y}t_H(m)
 \sum_{m_1 = 1}^{\infty} (\Lambda(m_1) - 1) e^{-m_1 / N}
   \sum_{m_2 = 1}^{\infty} e^{-m_2 / N}
   \begin{cases}
     1 & \text{if $m_1 + m_2 = m+N$} \\
     0 & \text{otherwise}
   \end{cases} \\
\notag
  &=
   \sum_{m= -H}^{y}t_H(m)
    e^{- (m +N) / N}
    \sum_{m_1 = 1}^{m+N- 1} (\Lambda(m_1) - 1)
\\
\label{cancellation}
  &=
   \sum_{n= N-H}^{N+y}
   e^{- n / N} t_H(n-N)
     (\psi(n - 1) - (n - 1)),
\end{align}
since the condition $m_1 + m_2 = m+N$ implies that both variables are $<m+N$.
Now $\psi(n) = \psi(n - 1) + \Lambda(n)$, so that
\begin{align*}
\sum_{n= N-H}^{N+y}e^{-n / N} t_H(n-N) (\psi(n - 1) - (n - 1))
  &=
\sum_{n= N-H}^{N+y}e^{-n / N} t_H(n-N) (\psi(n) - n)
  +
  \Odi{E},
\end{align*}
say, where, by the Brun-Titchmarsh theorem, 
we have $E \ll H^2 \log N$ if $0\le H+y \le N^\eps$ and $E \ll H^2 $
otherwise.
By \eqref{start-third-step}-\eqref{cancellation} and the previous
equation, we have
\begin{align*}
  \int_{-1/2}^{1/2}    T_H(N,y;-\alpha)\frac{\widetilde{R}(\alpha)}z \, \dx \alpha
  &=
\sum_{n= N-H}^{N+y}e^{-n / N} t_H(n-N) (\psi(n) - n)
  +
  \Odi{H^{3/2} N^{1/2} (\log N)^{1/2}}.
\end{align*}
Hence \eqref{I2-eval} is proved.
\end{Proof}

We need now the following lemma which is an extension of Lemma 2.6 of
\cite{LanguascoZ2012a}.
\begin{Lemma}
\label{sum-integral-modificato}
Let $M > 1$ be a real number. We have that
\[
\sum_{n = 1}^M (\psi(n)-n)
=
-
\sum_{\rho} \frac{M^{\rho + 1}}{\rho (\rho + 1)}
+
\Odi{M},
\]
where $\rho$ runs over the non-trivial zeros of the Riemann zeta function $\zeta(s)$.
\end{Lemma}
\begin{Proof} 
The case when $M>1$ is an integer was proved in   Lemma 2.6 of
\cite{LanguascoZ2012a}. 
Let $M > 1$ be a non-integral real number. Hence 
\begin{equation}
\label{first-compare}
\sum_{n = 1}^M (\psi(n)-n) 
=
\sum_{n = 1}^{\lfloor M \rfloor} (\psi(n)-n)
=
-
\sum_{\rho} \frac{{\lfloor M \rfloor}^{\rho + 1}}{\rho (\rho + 1)}
+
\Odi{M},
\end{equation}
by  Lemma 2.6 of \cite{LanguascoZ2012a}.
Writing $\rho=\beta+i \gamma$, we have
\[
 \frac{{\lfloor M \rfloor}^{\rho + 1}-{ M }^{\rho + 1}}{ \rho + 1}
\ll
M^{\beta +1}
\min\Bigl(\frac{1}{M}; \frac{1}{\vert \rho+1 \vert}\Bigr)
\]
and hence, by the zero-free region and the Riemann-von Mangoldt estimate, we obtain
\begin{align}
\label{compare}
\sum_{\rho} \frac{{\lfloor M \rfloor}^{\rho + 1}-{ M }^{\rho + 1}}{\rho (\rho + 1)}
\ll
\sum_{\vert \rho\vert \le M} \frac{M^{\beta }}{\vert \rho \vert}
+
\sum_{\vert \rho\vert > M} \frac{M^{\beta +1 }}{\vert \rho \vert^2}
 = 
 \odi{M}.
\end{align}
By \eqref{first-compare}-\eqref{compare}, Lemma \ref{sum-integral-modificato} follows.
\end{Proof}

\begin{Lemma}
\label{sum-integral-pesato}
Let $N$ be a large integer and $2\leq H \leq N$. We have that
\[
\sum_{n = N-H}^{N+H} t_H(n-N)(\psi(n)-n)
=
-
\sum_{\rho} \frac{(N+H)^{\rho +2} - 2 N^{\rho+2} +(N-H)^{\rho +2} }{\rho (\rho + 1)(\rho + 2)}
+
\Odi{HN},
\]
where $\rho$ runs over the non-trivial zeros of the Riemann zeta function $\zeta(s)$.
\end{Lemma} 
\begin{Proof}
A direct computation shows 
\begin{align}
\notag
\sum_{n = N-H}^{N+H} t_H(n-N)(\psi(n)-n)
&=
 \sum_{m = 0}^{H} t_H(m)(\psi(N+m)-(N+m))
\\
\notag
&\hskip1cm +
 \sum_{m = 0}^{H} t_H(m)(\psi(N-m)-(N-m))
 - H (\psi(N)-N)
 \\
 \label{inizio-termine-zeri}
 & =
 \sum_{m = 0}^{H} t_H(m)a(m) 
 - H (\psi(N)-N),
\end{align}
where we have implicitly defined $a(m) := \psi(N+m)+\psi(N-m)-2N$.
By partial summation we have
\begin{equation*}
 \sum_{m = 0}^{H} t_H(m)a(m)
=
-H a(0) +
\int_0^H \sum_{m=0}^t a(m)\ \dx t.
\end{equation*}
It is easy to see that $a(0)= 2 (\psi(N)-N)$ and  that
\begin{align*}
\sum_{m=0}^t a(m) &= 
\sum_{n = N-t}^{N+t} (\psi(n)-n)  + (\psi(N)-N) 
\\&
=
\sum_{n = 1}^{N+t} (\psi(n)-n) 
-
\sum_{n = 1}^{N-t} (\psi(n)-n)
+
 (\psi(N)-N)  
 +^\prime (\psi(N-t)-(N-t))  
\\&
=  
- 
\sum_{\rho} \frac{(N+t)^{\rho + 1}-(N-t)^{\rho + 1}}{\rho (\rho + 1)}
+
\Odi{N},
\end{align*}
where $+^\prime$ indicates that the term is present only when $t$ is
an integer and, in the last equality, we used Lemma \ref{sum-integral-modificato}
and the Prime Number Theorem.
Summing up, exploiting the absolute convergence of the series
over the zeros of the Riemann zeta function $\zeta(s)$, we obtain that
\begin{align}
\notag
 \sum_{m = 0}^{H} t_H(m)a(m)
&=
-
\sum_{\rho} \frac{1}{\rho (\rho + 1)}
\int_0^H 
\Bigl( (N+t)^{\rho + 1}-(N-t)^{\rho + 1} \Bigr)
\ \dx t
+
\Odi{HN}
\\
 \label{fine-termine-zeri}
 &
=
-
\sum_{\rho} \frac{(N+H)^{\rho +2} - 2 N^{\rho+2} +(N-H)^{\rho +2} }{\rho (\rho + 1)(\rho + 2)}
+
\Odi{HN}.
\end{align}
Inserting \eqref{fine-termine-zeri} in \eqref{inizio-termine-zeri} and using the Prime Number Theorem,
Lemma \ref{sum-integral-pesato} follows.
\end{Proof}

\section{Proof of Theorem \ref{Main-Th}}

We will get Theorem \ref{Main-Th} as a consequence of
Theorem \ref{average-Th}.
By partial summation we have
\begin{align}
\notag
\sum_{n=N-H}^{N+H}
& t_H(n-N)
\Bigl(
R(n) - (2\psi(n)-n)
\Bigr)
=
\sum_{n=N-H}^{N+H}
e^{n/N} t_H(n-N)
\Bigl\{
\Bigl(
R(n) - (2\psi(n)-n)
\Bigr)
e^{-n/N}
\Bigr\}
\\
&
\notag
=
e^{(N+H)/N}
\sum_{n=N-H}^{N+H}
e^{-n/N} t_H(n-N) 
\Bigl(
R(n) - (2\psi(n)-n)
\Bigr)
\\
&
\qquad
\label{psum}
-\frac{1}{N}
\int_{N-H}^{N+H}
\Bigl\{
\sum_{n=N-H}^{w}
e^{-n/N} t_H(n-N) 
\Bigl(
R(n) - (2\psi(n)-n)
\Bigr) 
\Bigr\}
e^{w/N}
\ \dx w
+\Odi{1}.
\end{align}
Inserting  \eqref{average}-\eqref{average-H} in \eqref{psum} we get
\begin{align*} 
\sum_{n=N-H}^{N+H}
t_H(n-N) 
\Bigl(
R(n) - (2\psi(n)-n)
\Bigr)
\ll
HN \Bigl(\log \frac{2N}{H}\Bigr)^2 + H^2 (\log N)^2 \log(2H) 
\end{align*}
and hence
\begin{align}
\notag
\sum_{n=N-H}^{N+H}
 t_H(n-N)  R(n)
&=
\sum_{n=N-H}^{N+H}
 t_H(n-N) n
+
2\sum_{n=N-H}^{N+H} t_H(n-N) 
(\psi(n)-n)
\\
\label{GY-improved}
&
+\Odig{
HN \Bigl(\log \frac{2N}{H}\Bigr)^2 + H^2 (\log N)^2 \log(2H) 
}.
\end{align}
A direct calculation proves that
\[
\sum_{n=N-H}^{N+H}  t_H(n-N) n =  H^2N
\]
and hence 
Theorem \ref{Main-Th} now follows inserting such an identity and  
Lemma \ref{sum-integral-pesato} in \eqref{GY-improved} and dividing by $H$.

\section{Proof of  Theorem \ref{average-Th}}
Assume  $N \geq 2$, $1\le H \le N$, $y \in [-H, H]$  and
let $\alpha\in[-1/2,1/2]$.
Writing
$\widetilde{R}(\alpha)=\widetilde{S}(\alpha) - 1/z$,
recalling \eqref{T-def}
we have
\begin{align}
\notag
&\sum_{n=N-H}^{N+y}
e^{-n/N}t_H(n-N)R(n)
=
\int_{-\frac{1}{2}}^{\frac{1}{2}}
\widetilde{S}(\alpha)^2
T_H(N,y;-\alpha)
\ \dx \alpha
\\
\notag
&
=
\int_{-\frac{1}{2}}^{\frac{1}{2}}
\frac{T_H(N,y;-\alpha)}{z^2}
\ \dx \alpha
+
2
\int_{-\frac{1}{2}}^{\frac{1}{2}}
\frac{T_H(N,y;-\alpha)\widetilde{R}(\alpha)}{z}
\ \dx \alpha
+
\int_{-\frac{1}{2}}^{\frac{1}{2}}
T_H(N,y;-\alpha)\widetilde{R}(\alpha)^{2}
\ \dx \alpha
\\
\label{circle}
&
= I_{1}(y) + I_{2}(y) +I_{3}(y),
\end{align}
say.

\paragraph{\textbf{Evaluation of $I_{1}(y)$}}
By Lemma \ref{residue}
we obtain
\begin{align}
\notag
 I_{1}(y)
&
=
\int_{-\frac{1}{2}}^{\frac{1}{2}}
\frac{T_H(N,y;-\alpha)}{z^2}
\ \dx \alpha
=
\sum_{n=N-H}^{N+y}t_H(n-N)
\int_{-\frac{1}{2}}^{\frac{1}{2}}
\frac{e(-n\alpha)}{z^2}
\ \dx \alpha
\\
\label{I1-eval}
&=
\sum_{n=N-H}^{N+y}t_H(n-N)
\Bigl(
n e^{-n/N} +\Odi{1}
 \bigr)
=
\sum_{n=N-H}^{N+y} e^{-n/N}  t_H(n-N)
n 
+
\Odi{H(H+y+1)}.
\end{align}

\paragraph{\textbf{Estimation of $I_{2}(y)$}}
By \eqref{I2-eval} of  Lemma \ref{I2-lemma}
we obtain
\begin{align}
\label{I2-final}
I_{2}(y)
=
2
\sum_{n=N-H}^{N+y}  e^{-n / N} t_H(n-N)  (\psi(n) - n)   +
  \Odi{H^{3/2} N^{1/2} (\log N)^{1/2}}.
\end{align}

\paragraph{\textbf{Estimation of $I_{3}(y)$; $y\in[-H,H)$}}

Using  Lemmas \ref{T-behaviour} and \ref{LP-Lemma-mod} we have that
\begin{align}
\notag
 I_{3}(y)
&
\ll
\int_{-\frac{1}{2}}^{\frac{1}{2}}
\vert T_H(N,y;-\alpha)\vert
\vert \widetilde{R}(\alpha)\vert^{2}
\ \dx \alpha
\\ \notag
&
\ll
H^2
\int_{-\frac{1}{H}}^{\frac{1}{H}}
\vert \widetilde{R}(\alpha)\vert^{2}
\ \dx \alpha
+
H
\int_{\frac{1}{H}}^{\frac{1}{2}}
\frac{\vert \widetilde{R}(\alpha)\vert^{2}}{\alpha}
\ \dx \alpha
+
H
\int_{-\frac{1}{2}}^{-\frac{1}{H}}
\frac{\vert \widetilde{R}(\alpha)\vert^{2}}{\vert \alpha \vert}
\ \dx \alpha
\\
\notag
&
\ll
H N \Bigl(\log \frac{2N}{H}\Bigr)^2 
+
H
\sum_{k=0}^{\Odi{\log 2H }}
\frac{H}{2^{k}}
\int_{\frac{2^{k}}{H}}^{\frac{2^{k+1}}{H}}
\vert \widetilde{R}(\alpha)\vert^{2}
\ \dx \alpha
\\
\notag&
\ll
H N \Bigl(\log \frac{2N}{H}\Bigr)^2 
+
H
\sum_{k=0}^{\Odi{\log 2H}}
\frac{H}{2^{k}} N \frac{ 2^{k+1}}{H}  \Bigl(\log \frac{2^{k+2}N}{H}\Bigr)^2 
\\
\notag&
\ll
H N \Bigl(\log \frac{2N}{H}\Bigr)^2 
+
H N (\log  N )^2
\sum_{k=0}^{\Odi{\log 2H}} 1 
\\
\label{I3-estim}
&
\ll
H N  (\log  N )^2  \log (2H)
.
\end{align}

\paragraph{\textbf{Estimation of $I_{3}(H)$}}

Using  Lemmas \ref{T-behaviour} and \ref{LP-Lemma-mod} we have that
\begin{align}
\notag
 I_{3}(H)
&
\ll
\int_{-\frac{1}{2}}^{\frac{1}{2}}
\vert T_H(N,H;-\alpha)\vert
\vert \widetilde{R}(\alpha)\vert^{2}
\ \dx \alpha
\\ \notag
&
\ll
H^2
\int_{-\frac{1}{H}}^{\frac{1}{H}}
\vert \widetilde{R}(\alpha)\vert^{2}
\ \dx \alpha
+ 
\int_{\frac{1}{H}}^{\frac{1}{2}}
\frac{\vert \widetilde{R}(\alpha)\vert^{2}}{\alpha^2}
\ \dx \alpha
+ 
\int_{-\frac{1}{2}}^{-\frac{1}{H}}
\frac{\vert \widetilde{R}(\alpha)\vert^{2}}{\alpha^2 }
\ \dx \alpha
\\
\notag
&
\ll
H N \Bigl(\log \frac{2N}{H}\Bigr)^2 
+
\sum_{k=0}^{\Odi{\log 2H }}
\frac{H^2}{4^{k}}
\int_{\frac{2^{k}}{H}}^{\frac{2^{k+1}}{H}}
\vert \widetilde{R}(\alpha)\vert^{2}
\ \dx \alpha
\\
\notag&
\ll
H N \Bigl(\log \frac{2N}{H}\Bigr)^2 
+ 
\sum_{k=0}^{\Odi{\log 2H}}
\frac{H^2}{4^{k}} N \frac{ 2^{k+1}}{H}  \Bigl(\log \frac{2^{k+2}N}{H}\Bigr)^2 
\\
\notag&
\ll
H N \Bigl(\log \frac{2N}{H}\Bigr)^2 
+ 
HN
\sum_{k=0}^{\Odi{\log 2H}}
\frac{1}{2^{k}}  \Bigl(k+1+ \log \frac{2N}{H}\Bigr)^2 
\\
\label{I3-estim-H}
&
\ll
H N \Bigl(\log \frac{2N}{H}\Bigr)^2 
.
\end{align}

\paragraph{\textbf{End of the proof}}

Inserting \eqref{I1-eval}-\eqref{I3-estim}
into \eqref{circle}, for every $y\in [-H,H)$,  we immediately have
\begin{align*}
\sum_{n=N-H}^{N+y}
& e^{-n/N} t_H(n-N) R(n)
=
\sum_{n=N-H}^{N+y}
e^{-n/N}  t_H(n-N) n 
\\&
+
2
\sum_{n=N-H}^{N+y}
e^{-n/N} t_H(n-N)  (\psi(n)-n)
+
\Odig{
H N (\log N )^2   \log (2H)
}.
\end{align*}
Hence
\[
\sum_{n=N-H}^{N+y}
e^{-n/N} t_H(n-N) 
\Bigl(
R(n) -(2\psi(n)-n)
\Bigr)
\ll
H N (\log N )^2  \log (2H)
\]
for every $y\in [-H,H)$.
Thus we can write
\begin{equation*}
\max_{y\in  [-H,H)}
\Bigl  \vert
\sum_{n=N-H}^{N+y}
e^{-n/N}  t_H(n-N) 
\Bigl(
R(n) -(2\psi(n)-n)
\Bigr)
\Bigr \vert
\ll
H N  (\log N )^2   \log (2H).
\end{equation*}
The $y=H$ case follows analogously using \eqref{I3-estim-H} instead of \eqref{I3-estim}.
Dividing by $H$, Theorem \ref{average-Th} is proved.

 \section{About the order of magnitude of the zero-depending term}
 \label{S-order}

Let us define
\begin{align*}
S &:= 
\sum_{\rho} \frac{(N+H)^{\rho +2} - 2 N^{\rho+2} +(N-H)^{\rho +2} }
{\rho (\rho + 1)(\rho + 2)}.
\end{align*}
%
%
By Lemma \ref{sum-integral-pesato} we have
\begin{align*} 
 S
= 
- \sum_{n = N-H}^{N+H} t_H(n-N)(\psi(n)-n)
+
\Odi{HN}.
\end{align*}
Assuming RH, we have $\psi(n)-n \ll n^{1/2} (\log n)^2$
and hence
\begin{align*} 
S
&\ll
H (\log N)^2 \sum_{n = N-H}^{N+H} n^{1/2}
+ 
HN
\ll
H (\log N)^2 \Bigl((N+H)^{3/2} - (N-H)^{3/2} \Bigr) 
+ 
HN
\\&
\ll
H^2 N^{1/2} (\log N)^2 + HN.
\end{align*}
Dividing by $H$,
the expected order 
of magnitude of the  the second difference term in Theorem \ref{Main-Th}
is, under the assumption of RH,  $\ll HN^{1/2} (\log N)^2 +N$. 
%
%
The same strategy works in the unconditional case
too. The Prime Number Theorem in the form
$\psi(n)-n \ll n \exp(-c (\log n)^{3/5}(\log \log n)^{-1/5})$,
where $c>0$ is an absolute constant, leads to the final estimate 
\begin{align*}
S
&\ll
H \exp(-c_1 (\log N)^{3/5}(\log \log n)^{-1/5}) \sum_{n = N-H}^{N+H} n
+HN
\\
&\ll
H^2 N \exp(-c_1 (\log N)^{3/5}(\log \log n)^{-1/5})
+ HN,
\end{align*}
where $c_1>0$ is an absolute constant.
Dividing by $H$,
the expected order 
of magnitude of the  the second difference term in Theorem \ref{Main-Th} is $\ll HN \exp(-c_1 (\log N)^{3/5}(\log \log n)^{-1/5})+N$. 

%

\vskip 1cm
\noindent
Alessandro Languasco, Dipartimento di Matematica  ``Tullio Levi-Civita'', Universit\`a
di Padova, Via Trieste 63, 35121 Padova, Italy; languasco@math.unipd.it

\medskip
\noindent
Alessandro Zaccagnini, Dipartimento di Matematica e Informatica, Universit\`a di Parma, Parco
Area delle Scienze 53/a, 43124 Parma, Italy;
alessandro.zaccagnini@unipr.it

\end{document}